\documentclass[11pt,psamsfonts,twoside]{article}
\usepackage{mathrsfs}
\usepackage{bbm}
\usepackage{amsmath}
\usepackage{amsmath,amsthm}
\usepackage{amssymb,amscd}
\usepackage{amsfonts,amsbsy}
\usepackage{fancyhdr,graphicx}
\usepackage[dvips]{psfrag}
\usepackage{indentfirst}
\usepackage{xcolor}
\usepackage{setspace}

\numberwithin{equation}{section}
\newtheorem {proposition} {Proposition}[section]
\newtheorem {theorem}     [proposition]{Theorem}
\newtheorem {corollary}   [proposition]{Corollary}
\newtheorem {lemma}       [proposition]{Lemma}

\newtheorem {remark}      [proposition]{Remark}

\newtheorem {definition}  [proposition]{Definition}
\allowdisplaybreaks

\newcommand{\supp}{{\rm supp}}
\newcommand{\diva}{{\rm div}}

\begin{document}
\setlength{\parindent}{4ex} \setlength{\parskip}{1ex}
\setlength{\oddsidemargin}{12mm} \setlength{\evensidemargin}{9mm}
\title{\textsc{\text{ILL-POSEDNESS OF THE INCOMPRESSIBLE} \\ NAVIER-STOKES EQUATIONS IN $\dot{F}^{-1,q}_{\infty}(\mathbb{R}^3)$}}
\author{\textsc{Chao Deng and Xiaohua Yao}}

\date{}
\maketitle
\begin{abstract}
In this paper, authors show the ill-posedness of 3D incompressible
Navier-Stokes equations in the critical Triebel-Lizorkin spaces
$\dot{F}^{-1,q}_{\infty}(\mathbb{R}^3)$  for any $q>2$ in the sense
 that arbitrarily small initial data of
$\dot{F}^{-1,q}_{\infty}(\mathbb{R}^3)$ can lead the corresponding
solution to become arbitrarily large after an arbitrarily short
time. Thus extends Bourgain and Pavlovi$\acute{c}$'s work
\cite{BP08}. In view of the well-posedness of 3D-incompressible
Navier-Stokes equations in $BMO^{-1}$  ( i.e. the Triebel-Lizorkin
space $ \dot{F}^{-1,2}_{\infty}(\mathbb{R}^3)$ ) by Koch and Tataru,
our work completes a  {\it dichotomy} of well-posedness and
ill-posedness in the Triebel-Lizorkin space framework depending on
$q=2$ or $q > 2$.
 \end{abstract}
 {\small \noindent{\bf Keywords:}\,  Navier-Stokes equations; Triebel-Lizorkin space; Well-posedness; Ill-posedness.}
\maketitle
\section{Introduction}
 In this article, we are concerned with the following incompressible 3D Navier-Stokes equations (NS):
   \begin{align}\label{NS:1.1}
  \left\{\begin{aligned}
   &\partial_t u + \Delta u + u\cdot \nabla u +\nabla
   p=0,\\
   &\diva \,u=0,\\
   &u(x,0)=u_0(x),
     \end{aligned}\right.
 \end{align}
 where $(x,t)\in\mathbb{R}^3\times(0,\infty)$, ${u}(x,t)=\!(u^1(x,t), u^2(x,t),\!u^3(x,t))$
 are unknown vector functions, $p(x,t)$ is unknown scaler function,  and ${u}_0(x)$ is a given vector function satisfying
 divergence free condition $\nabla\cdot{u}_0=0$.

 Mathematical study on the existence and uniqueness of the incompressible Navier-Stokes equations has a long history.  In 1934, Leray \cite{LER34} first proved existence
 of global weak solution associated with any $L^2(\mathbb{R}^3)$ initial data by some weak compactness arguments.
 Until now, whether such a weak solution is unique and smooth or not is still a great open problem, See Fefferman \cite{Feff}, also Caffarelli, Kohn and Nirenberg \cite{CKN} and F.-H. Lin \cite{Lin} for partial regularity of suitable weak solution.  Beginning with a different method by semigroup and Picard's iteration,
 Fujita-Kato \cite{FUJK64} in 1964 established the local well-posedness of the N-S in ${H}^s(\mathbb{R}^3)$ for any $s\geq\frac{3}{2}-1$, and also global well-posedness for any small initial data of ${H}^s(\mathbb{R}^3)$. This
 remarkable approach can be adapted to  other various function spaces of initial data, see \cite{CAN04,Kato:1984471,LEM02,Meyer99} for expositions and references therein. In particular, an
interesting result that should be mentioned is due to Koch-Tataru
\cite{KOCT01}, where they proved that the N-S equation is well-posed
in $BMO^{-1}$ ( i.e. the Triebel-Lizorkin space $
\dot{F}^{-1,2}_{\infty}(\mathbb{R}^3)$ ).

 On the other hand, recently, Bourgain and Pavlovi$\acute{c}$ in \cite{BP08} proved the ill-posedness of the incompressible 3D-Navier-Stokes equations in the largest scaling invariant
Besov function space $\dot{B}_\infty^{-1,\infty}(\mathbb{R}^3)$ in
the sense that arbitrarily small initial data in the space
$\dot{B}_\infty^{-1,\infty}(\mathbb{R}^3)$
 can lead the solution to become arbitrarily large after an arbitrarily short time. Notice that $BMO^{-1}\varsubsetneq \dot{B}_\infty^{-1,\infty}(\mathbb{R}^3)$,
 clearly, there exists a slight difference between the well-posed space and the ill-posed space. Therefore, it would be an interesting problem whether there exist
 strictly smaller ill-posed spaces than $\dot{B}_\infty^{-1,\infty}(\mathbb{R}^3)$.   Such an improvement of illposedness was obtained  by
 Yoneda \cite{Yoneda:2010} in a logarithmic type Besov space near $BMO^{-1}$.  Motivated by that  $BMO^{-1}$ is identical to the end-pointed  Triebel-Lizorkin space
 $\dot{F}_\infty^{-1, 2}(\mathbb{R}^3)$, in the current paper, the authors  further study the interesting problem in a critical Triebel-Lizorkin space
 framework and  show that the the NS eqaution \eqref{NS:1.1} is ill-posed in $\dot{F}^{-1,q}_{\infty}(\mathbb{R}^3)$ for any $q>2$, which are all strictly
 smaller than $\dot{B}_\infty^{-1,\infty}(\mathbb{R}^3)$ ( See Theorem \ref{Thm:1.3} below). Hence our work along with \cite{KOCT01} establishes a dichotomy of well-posedness and ill-posedness in the Triebel-Lizorkin space framework depending on $q=2$ or $q > 2$.

  For the end, let us first recall the definitions of homogeneous Besov spaces/Triebel -Lizorkin spaces. Let $\varphi(\xi)=\varphi(|\xi|)\ge0$ be a real-valued
  smooth function such that
  \begin{align}\label{NS:1.2}
  \left\{\begin{aligned}
   &\supp\,\varphi\subset\!\{\xi\in\mathbb{R}^3;\;
   {5}/{8}\le|\xi|\le {7}/{4}\},\\
   &\varphi\equiv1 \text{ in } \{\xi\in\mathbb{R}^3;\;  {7}/{8}\le|\xi|\le
    {5}/{4}\}, \\
   &\sum{\!}_{j}\;\varphi(2^{-j}\xi)=1  \quad\text{ for }
   \xi\in\mathbb{R}^3\backslash\{0\}.
   \end{aligned}
   \right.
  \end{align}
For any tempered distribution $f$ and $i,j\in\mathbb{Z}$, define the dyadic block as follows:
 \begin{align}\label{NS:1.3}
 \Delta_{j}f(x)=\varphi(2^{-j}\nabla)f(x)\quad \text{ and\;\; }
  \Delta_{i}\Delta_jf\equiv0 \text{ if } |i-j|\ge2.
\end{align}
 In order to exclude nonzero polynomials in homogeneous Besov spaces and Triebel-Lizorkin spaces,
 it is natural to use $Z'(\mathbb{R}^3)$ to denote the subset of tempered
 distribution $f\in {S}\hskip.01cm '(\mathbb{R}^3)$ modulo all polynomials set $P(\mathbb{R}^3)$, i.e.
 $Z'(\mathbb{R}^3)=S\hskip.01cm '(\mathbb{R}^3)/P(\mathbb{R}^3).$

 Now we are ready to  give the definitions of Triebel-Lizorkin spaces $\dot{F}^{-1,q}_{\infty}(\mathbb{R}^3)$ and Besov space $\dot{B}^{-1,\infty}_{\infty}(\mathbb{R}^3)$,
 also see
\cite{Triebel:1978} for a detailed exposition about other general
spaces $\dot{F}^{s,q}_{p}(\mathbb{R}^3)$ and
$\dot{B}^{s,q}_{p}(\mathbb{R}^3)$.

\begin{definition}\label{Def:1.1}
 For $1<q<\infty$, we define $\dot{F}^{-1,q}_{\infty}(\mathbb{R}^3)$
 as the following set so that
 \begin{align}\label{NS:1.4}
   \dot{F}^{-1, q}_{\infty}(\mathbb{R}^3)
     =
   \Big\{
   & \,f\, |
     \; f\!\in\! {Z}'(\mathbb{R}^3),\, \exists\, \{f_k(x)\}_{k\in\mathbb{Z}}\; {\rm\; s.t. \;} f\! =
     \!\sum{_{k\in\mathbb{Z}}}\,
  \Delta_k f_{k}(\cdot)\;\; \nonumber\\
   &\text{in
   } Z'(\mathbb{R}^3) \;\;\text{and }
     \Big\|  \|\{2^{-k}|f_k(\cdot)|\} \|_{\ell^{\,q}} \Big\|_{L^\infty(\mathbb{R}^3)}<\infty
   \Big\}
 \end{align}
 and the corresponding norm is defined by
 $$\|f \|_{\dot{F}^{-1,q}_{\infty}(\mathbb{R}^3)} = \inf \| \{2^{-k}|f_k(\cdot)|\}_{\ell^{\,q}}
 \|_{L^\infty(\mathbb{R}^3)},$$
 where the infimum is taken over all admissible representations in the sense of \eqref{NS:1.4}.
 Meanwhile, we denote by
 $\dot{B}^{-1,\infty}_{\infty}(\mathbb{R}^3)$ the set of distribution
 $f\!\in\!Z{'}(\mathbb{R}^3)$ such  that
 \begin{align}\label{NS:1.5}
   \|f\|_{\dot{B}^{-1,\infty}_{\infty}(\mathbb{R}^3)}=\sup_{t>0}\sqrt{t}\,\|e^{t\Delta}f \|_{L^\infty(\mathbb{R}^3)}<\infty.
 \end{align}
 \end{definition}

 \begin{remark}
 It is known that $\dot{F}^{-1,2}_{\infty}(\mathbb{R}^3)=BMO^{-1}(\mathbb{R}^3)$ and  $BMO^{-1}(\mathbb{R}^3)$ has
 the following equivalent Carleson measure characterization {\rm (cf. \cite{KOCT01})}:
 \begin{align}\label{NS:1.5B}
   \big\|\,f\big\|_{BMO^{-1}(\mathbb{R}^3)}=
   \sup_{x\in\mathbb{R}^3,R>0}\Big(\frac{1}{|B_R(x)|}\int_{0}^{R^2}\!\!\!\int_{B_R(x)}|e^{t\Delta}f(y)|^2dtdy\Big)^{\frac{1}{2}}<\infty.
 \end{align}
 Moreover,  we remark that for any $1<q<\infty$,
 $$\dot{F}^{-1,q}_{\infty}(\mathbb{R}^3)\hookrightarrow
 \dot{B}^{-1,\infty}_{\infty}(\mathbb{R}^3).$$
 \end{remark}

As usual, we first write \eqref{NS:1.1} into
 the following equivalent integral equations:
\begin{align}
  &u=e^{t\Delta}u_0-B(u,u),\label{DLY-2}
  \end{align}
  where $\mathbb{P}$ denotes the Leray projection operator
  $Id-\nabla\frac{1}{\Delta}\diva$,  and the bilinear term  $B(u,v)$ is defined by
\begin{align}
  \label{DLY-3}
 B(u,v):=\int_{0}^te^{(t-\tau)\Delta}\mathbb{P}(u\cdot\nabla v) d\tau.
 \end{align}

For any $u\in L^2_{loc}(\mathbb{R}^3\times\mathbb{R})$, we define that
 \begin{align}\label{NS:1.6}
   \|u\|_{X_T}:=\sup_{x_0\in\mathbb{R}^3, 0<R^2<T}
   \Big(\frac{1}{|B_R(x_0)|}\int_0^{R^2}\int_{B_R(x_0)}|u|^2dtdx\Big)^{\frac{1}{2}},
 \end{align}
 and\vspace*{-2ex}
 \begin{align}\label{NS:1.7}
   \|u\|_{\mathscr{E}_T}:= \sup_{0<t<T}t^{\frac{1}{2}}\|u\|_{L^\infty(\mathbb{R}^3)} + \|u\|_{X_T}.
 \end{align}

  Recall that $B(u,v)$
 satisfies the following {\it a-priori} bilinear estimates:
  \begin{align}\label{NS:1.8}
  \|B(u,v)\|_{\mathscr{E}_T}\lesssim \|u\|_{\mathscr{E}_T}\|v\|_{\mathscr{E}_T}.
  \end{align}
 Applying boundedness property of $\mathbb{P}$ in $\dot{F}^{-1,2}_\infty(\mathbb{R}^3)$ and
 decay estimate for heat kernel, following similar argument as in
 \cite[Lemma 16.3]{LEM02}, we have
 \begin{align}\label{NS:1.9}
  \|B(u,v)\|_{L^\infty_T\dot{F}^{-1,2}_\infty}\lesssim \|u\|_{\mathscr{E}_T}\|v\|_{\mathscr{E}_T}.
  \end{align}
Based on \eqref{NS:1.8}, Koch and Tataru \cite{KOCT01} established
the existence of solutions to the Navier-Stokes equation in
$BMO^{-1}(\mathbb{R}^3)$. By using \eqref{NS:1.8} and
\eqref{NS:1.9}, we will further prove the following ill-posedness in
$\dot{F}^{-1,q}_{\infty}(\mathbb{R}^3)$ for any $q>2$.

\begin{theorem}\label{Thm:1.3}
 For any $q>2$ and $\delta>0$,  there exists a solution $u$ to  system \eqref{NS:1.1} with some initial data $u_0\in
\dot{F}^{-1,q}_{\infty}(\mathbb{R}^3)$ satisfying
\begin{align*}
  \|u_0\|_{\dot{F}^{-1,q}_{\infty}(\mathbb{R}^3)}\lesssim\delta
  \end{align*}
and $\diva\, u_0=0$ such that for some $0<T<\delta$,
\begin{align*}
  \|u(T)\|_{\dot{F}^{-1,q}_{\infty}(\mathbb{R}^3)}\gtrsim\frac{1}{\delta}.
\end{align*}
\end{theorem}



 { This paper is organized as follows:}
 In Section 2, we first construct a very special initial data and list some
 necessary remarks and lemmas. In section 3, we establish all the desired estimates
 about the first and second approximation terms which will be used in controlling the remainder term.
 Finally, combining all the a-priori estimates we prove ill-posedness of the Navier-Stokes equations.

  \medskip
 \noindent\textit{Notations}:
 Throughout this paper, we shall use $C$ and $c$ to denote generic constants
  and may change from line to line. Both $\mathcal{F}f$ and $\widehat{f}$  stand for Fourier
 transform of $f$ with respect to space variable, while
 $\mathcal{F}^{-1}$ stands for the inverse Fourier transform.  We
 denote $A\le{CB}$ by $A\lesssim B$ and $A\lesssim{B}\lesssim{A}$
  by $A\sim{B}$. For any $1\le p\le\infty$, we denote
   $L^q(0,T)$,  $L^q(0,\infty)$ and $L^p(\mathbb{R}^3)$ by
   $L^q_T$,  $L^q_t$ and $L^p_x$, respectively.
  Later on, we also use $\dot{F}^{s,q}_p$ to denote $\dot{F}^{s,q}_p(\mathbb{R}^3)$
  if there is no confusion about the domain, and similar conventions
  are applied.  {For simplicity, we denote by $B_R(x)$ the ball centered at $x$ of radius $R$.}

\section{Construction of initial data}
 For any $\delta>0$, we define the initial data as follows:
 \begin{align}\label{NS:2.1}
  u_0(x)
  \!=\!\frac{Q}{\sqrt{r}}
     \sum_{s=1}^{r} \!\Big(
     \cos(k_sx){\Psi}_{1}  \!- \!|k_s|\sin(k_sx){\Psi}_{2}   \!+\!  \cos(k_s'x){\Psi}_{3}   \!+\!  |k_s|\sin(k_s'x){\Psi}_4
     \Big),
 \end{align}
 where
  \begin{align}\tag{H}
  \left.\begin{aligned}
 &k_s\!=\!(0,\,2^{\frac{(s+1)(s+2m_0)}{2}},\,0),\;\;
  k_s'\!=\!(2^3,\,-2^{\frac{(s+1)(s+2m_0)}{2}},\,0)\;\; \text{ with  $s\!=0, 1,2,\cdots, r$,} \medskip\\
  &{\Psi}_{1} \!=\!(0,-\partial_3\psi,\partial_2\psi),\!\!\!,\;
  {\Psi}_{2}\!=\!(0,0,\psi),\!\!\!,\;
   {\Psi}_{3} \!=\!(\partial_2\psi,-\partial_1\psi,0),\!\!\!,\;
  {\Psi}_4\!=\!(\psi,\frac{2^3\psi}{|k_s|},0)\!\!\!,\\
  &\psi(x) \text{ satisfies:  $\widehat{\psi}(\xi)\!=\!\widehat{\psi}(|\xi|)\ge0,\;
    \ \textrm{supp}\,\widehat{\psi}\subset
     B_{\frac{1}{4}}(0)$ and
     $\|\widehat{\psi}\|_{L^1_\xi}=1$}\end{aligned}\right\}
 \end{align}
 and $Q$, $r$, $m_0$ {will be chosen sufficiently large} {according to the size of} $\delta$.
%
\vskip0.5cm
 \begin{remark}\label{Rem:2.1} {\rm
    From the above assumptions we have the following observations:
  \begin{itemize}
  \vspace*{-1ex}
  \item[(i)]
  From \eqref{NS:2.1} and (H), it is easy to check that
   $u_0$ is {\it smooth, real-valued, divergence free and $L^2_x$ finite}.
   \item[(ii)] For $\forall\, g\!\in\! S(\mathbb{R}^3)$ and  $\forall\, k\!\in\!\mathbb{Z}^3$ we have
   $\widehat{g}(\xi\!-\!k)=\mathcal{F}( e^{ikx}g(x))$, 
%
 which shows that
 \begin{align}\label{NS:2.2B}
   &\supp\,\mathcal{F}(\,(\,\cos(kx)\pm\sin(kx)\,)\psi(x)\,)\subset\! {B}_{\frac{1}{4}}(k)\cup B_{\frac{1}{4}}(-k).
  \end{align}
\item[(iii)] Recall that $|k_s|=2^{\frac{(s+1)(s+2m_0)}{2}}$, for any $s\in [0,r]\cap \mathbb{N}$, we
 denote $j_s=\log_2|k_s|$. Then for
 large enough $m_0$ ($m_0\ge\!7$), 
 %
  \begin{align}\label{NS:2.3}
      &   (\,
          \cos(k_sx)\psi,\,\cos(k_s'x)\psi
          \,)
        \!=
          (\,
         \Delta_{j_s} (\cos(k_sx)\psi), \Delta_{j_s}
         (\cos(k_s'x)\psi)
          \,).
  \end{align}
Similar arguments work well for $\sin(k_sx)\psi$ and $\sin(k_s'x)\psi$. 
  \item[(iv)]\;  From \eqref{NS:2.1} and (iii), we
   denote
    \begin{align} \label{NS:2.4}
   \hskip-.2cm\!f_\ell \!= \!\left\{
  \begin{aligned}
   &\!\!\cos(k_sx){\Psi}_{1}  \!-\! |k_s|\sin(k_sx){\Psi}_{2}   \!+\!  \cos(k_s'x){\Psi}_{3}   \!+\!  |k_s|\sin(k_s'x){\Psi}_4   \textit{\,if\,\,}
       \ell\!=\!j_s,\\
   &\!0 {\rm \;\;\;otherwise, }
  \end{aligned}\right.
 \end{align}
 then $u_0\!=\!\frac{Q}{\sqrt{r}} \displaystyle{\sum_{\ell\,\in \mathbb{Z}}} \Delta_\ell f_\ell$ and
 $\Big\{\frac{Qf_\ell}{\sqrt{r}}\Big\}_{\!\ell}$ is a decomposition of
 $u_0$ in the sense of \eqref{NS:1.4}.
 \end{itemize}
 }
 \end{remark}

 \begin{lemma}\label{NS:lem2.2}
 Let $a(D)$ be a 3-dimensional Fourier multiplier operator corresponding to the homogeneous symbol $a(\xi)$ of degree $m\ge 0$. Then  there exists some $c>0$ such that for
  any $j\in\mathbb{Z}$ and $t\ge0$, the following point-wise estimate holds
 \begin{align}\label{NS:2.5}
    \big|(a(D)e^{t\Delta} \Delta_j f)(x)\big|\lesssim 2^{mj}e^{-ct2^{2j}}
   (Mf)(x), \;\;\; 
   \end{align}
where $\Delta_{j}\!$ is defined in \eqref{NS:1.3} and $Mf$ is the
Hardy-Littlewood maximal function of $f$. 
 \end{lemma}
 \begin{proof}
 Let $K_j(t,x)$ be the kernel of  $a(D)e^{t\Delta} \Delta_j$.  Then by scaling we have
 $$K_j(t,x)=2^{(m+3)j}\int_{\mathbb{R}^3}e^{i 2^j x\xi}e^{-t2^{2j}|\xi|^2}a(\xi)\varphi(\xi)d\xi, $$
 which immediately yields $|K_j(t,x)|\lesssim 2^{(m+3)j}e^{-ct2^{2j}}.$  By integrating by
 parts we get
 \begin{align*}
     |2^j x|^{4} \big| K_j(t,x)\big |
     &\lesssim 2^{(3+m)j}\!\sum_{|\mu|=4} \Big| \int_{\mathbb{R}^3} e^{i\xi2^jx} \partial^\mu_\xi \big(e^{-t2^{2 j}|\xi|^{2}} a(\xi)\varphi(\xi)\big)\ d\xi \Big| \nonumber\\
    &\lesssim 2^{(3+m)j}\!\int_{\supp\varphi}\!
     e^{-t2^{2j}|\xi|^{2}}\big(t^{4}2^{8j}\!+\!1\,\big) d\xi\\
     &\lesssim 2^{(3+m)j} e^{-ct2^{2j}},
  \end{align*}
   where we used
  $e^{- t2^{2 j}|\xi|^{2}}t^{N}2^{2\alpha Nj}
 \le C_N e^{-ct2^{2 j}}$ for $\xi\in\supp\,\varphi$. Hence we obtain
 that
 $$|K_j(t,x)|\lesssim 2^{mj}e^{-ct2^{2 j}}2^{3j}(1+2^j|x|)^{-4},$$
 which concludes the desired estimate \eqref{NS:2.5}.  
 \end{proof}

 \begin{lemma}\label{NS:lem2.3}
  Let $\psi$ be defined in {\rm (H)}, $M\psi$ be the Hardy-Littewood maximal function of $\psi$ and $\theta\!=\!(M\psi)^2$. For any
  $(k,h)\!\in \mathbb{Z}^3\times\mathbb{Z}^3$ and $\min\{|h|,|k|,|h\!+\!k|\}\ge 6$, we
  denote 
  \begin{align}\label{NS:2.6}
    F(t,\tau,x;k,h) =
    a_1(D)\,e^{(t-\tau)\Delta}\big(a_2(D)e^{\tau\Delta}\psi_{k}\,a_3(D)e^{\tau\Delta}\psi_h\big)\,(x)\,\;\;\;\;
  \end{align}
  where $a_\ell(D)$ are  Fourier multipliers with homogeneous symbols $a_\ell(\xi)$ of degree $m_\ell\ge0$ $($$\ell=1,2,3$$)$, $\psi_k$
  is either $\cos(kx)\psi(x)$
  or $\sin(kx)\psi(x)$ and $\psi_h$ is either $\cos(hx)\psi(x)$ or
  $\sin(hx)\psi(x)$.
  Then there exist positive constants $c$ and $C$ such that 
   \begin{align}
       \label{NS:2.7}
  & |a_1(D)e^{t\Delta}\psi_{h}(x)|\le
  C|h|^{m_1}\,e^{-ct|h|^2}(M\psi)(x)\,\;\;\;\;\;\;
          \end{align}
          and
          \begin{align}
     \label{NS:2.8}
     &|\,F(t,\tau,x;k,h)\,|\le C \frac{|k\!+\!h|^{m_1}|k|^{m_2}|h|^{m_3}}{e^{c(t-\tau)|k+h|^2+c\tau(|k|^2+|h|^2)}}\,(M\theta)(x).
  \end{align}
  \end{lemma}
   \begin{proof}
  We first recall that $\psi$ satisfies: $\widehat{\psi}(\xi)\ge0$, $\textrm{supp}\,\widehat{\psi}\subset
     B_{\frac{1}{4}}(0)$ and $\|\widehat{\psi}\|_{L^1_\xi}=1$.
 %
 For any $h\in\mathbb{Z}^3$ and $|h|> 2^{2}$, we
 have\footnote{$[\log_2|h|]$ stands for the integer part of
$\log_2|h|$.} $[\log|h|]\ge 2$. Furthermore, we get
 $$\frac{7}{8}2^{[\log_2|h|]}+2^{-1}<|h|<\frac{5}{4}2^{[\log_2|h|]+1}-2.$$
 Similar to Remark \ref{Rem:2.1} (iii) and (iv), we get
  \begin{align}\label{NS:2.11}
  a_1(D)e^{t\Delta}\psi_h(x)= a_1(D)e^{t\Delta}(\Delta_{[\log_2|h|]}+\Delta_{[\log_2|h|]+1})\psi_h(x).
  \end{align}
 Therefore,  \eqref{NS:2.7} follows immediately from
\eqref{NS:2.5}, \eqref{NS:2.11} and $|\psi_h(x)|\!\le\!|\psi(x)|$,
i.e.
 \begin{align*}
 |a_1(D)e^{t\Delta}\psi_{h}(x)|
 &\lesssim \sum_{\sigma=0,1}| a_1(D)e^{t\Delta}\Delta_{[\log_2|h|]+\sigma} \psi_h(x)|\\
 &\lesssim |h|^{m_1}e^{-ct|h|^2}(M\psi)(x).
 \end{align*}
 Again by applying Lemma \ref{NS:lem2.2} repeatedly to
  $f=a_2(D)e^{\tau\Delta}\psi_{k}$, $f=a_3(D)e^{\tau\Delta}\psi_h$ and $f=a_2(D)e^{\tau\Delta}\psi_{k}\,a_3(D)e^{\tau\Delta}\psi_h$,
   we
 can prove the desired estimate \eqref{NS:2.8}.
%
%
 \end{proof}

%

\vskip0.3cm

\begin{lemma}\label{NS:lem2.4} Let $\mu\ge 0$ and $\ell=0,1$. We have the
following estimates
\begin{align}\label{NS:2.12}
   \sum_{s=1}^r e^{-ct|k_s|^2}t^\mu|k_s|^{2\mu}|k_{s-\ell}|\lesssim \Big( \sum_{s=1}^r e^{-\frac{c}{2}t|k_s|^2}|k_{s-\ell}|^2
    \Big)^{\frac{1}{2}}.
  \end{align}
 \end{lemma}
 \begin{proof}
 Noticing that $e^{-\frac{c}{2}t|k_s|^2}t^\mu|k_s|^{2\mu}\lesssim 1$ and squaring the left side of \eqref{NS:2.12} we
 get
\begin{align*}
    \big(\sum_{s=1}^r e^{-\frac{c}{2}t|k_s|^2}|k_{s-\ell}|\big)^2
  & \lesssim \sum_{s=1}^r e^{-\frac{c}{2}t|k_s|^2}(|k_{s-\ell}|^2 + |k_{s-\ell}|(|k_{s-1-\ell}|+\cdots+|k_{1-\ell}|))\\
&\lesssim \sum_{s=1}^r e^{-\frac{c}{2}t|k_s|^2}|k_{s-\ell}|^2,
\end{align*}
which concludes the desired estimates.
   \end{proof}


\section{ Analysis of ill-posedness }
 In this section, we will prove ``norm inflation" of the NS equaiton
 in $\dot{F}^{-1,q}_{\infty}$ with
 $q>2$.
 Following the
 ideas in \cite{BP08}, we rewrite
 the solution to the NS equations as a  summation of the first approximation terms, the second
 approximation terms and remainder terms, i.e.
\begin{align}\label{NS:3.1}
 &u=u_1 \!- u_{2} + y,
 \end{align}
 where
 $u_1=e^{t\Delta}u_0$ and $u_{2}\!=B(u_1,u_1).$
 Moreover, the remainder terms
 satisfy the following integral equations:
 \begin{align}\label{NS:3.2}
 y &= G_0+ G_1 - G_2 ,
 \end{align}
       on $(0,\infty)$ with the
initial conditions $y(0)=0$, $G_0  = B(u_{2},
u_1\!-\!u_2)+B(u_1,u_{2})$ and
\begin{align*}
  \,
  G_1  = B(y,u_2\!-\!u_1)+B(u_2\!-\!u_1,y), \,\;\;\;
  G_2  = B(y,y).
      \end{align*}

 In the following, we will establish the {\it
 a-priori} estimates for $u_0$, $u_{1}$, $u_2$ and $y$. Precisely,
 In Subsection 3.1 we estimate the small upper bounds of $u_0$ and $u_1$;
 In Subsection 3.2, we prove both upper bound and lower bound of $u_2$;
 In Subsection 3.3, we prove the upper bound of $y$; In
 Subsection 3.4, we complete the proof of Theorem \ref{Thm:1.3}.

\subsection{Estimates for initial data and the first approximation terms}
  In this subsection, we will estimate $u_0$ and $e^{t\Delta}u_0$.
 \smallskip
 \begin{lemma}\label{NS:lem3.1}\!
 For any initial $u_0$ defined in \eqref{NS:2.1} and any $q\!\ge 2$, we
 obtain that
  \begin{align}\label{NS:3.4}
   &{\|u_0\|_{\dot{F}^{-1,q}_{\infty}}+\|e^{t\Delta}u_0\|_{\dot{F}^{-1,q}_{\infty}}\!\lesssim\!
    Q\,r^{\frac{1}{q}-\frac{1}{2}}},
\end{align}
for some absolute  constant $c>0$.
 \end{lemma}
\begin{proof}
 In view of the construction of $u_0$ and \eqref{NS:2.4},
 we get
 $$u_0=\frac{Q}{\sqrt{r}}\sum_{\ell\in\mathbb{Z}}\Delta_\ell
 f_\ell,\;\quad
 e^{t\Delta}u_0=\frac{Q}{\sqrt{r}}\sum_{\ell\in\mathbb{Z}}\Delta_\ell\, e^{t\Delta}
 f_\ell.$$
 By Definition \ref{Def:1.1} and \eqref{NS:2.4}, we have
  \begin{align}
  \begin{aligned}
   \|u_0\|_{\dot{F}^{-1,q}_{\infty}}
    &\lesssim\frac{Q}{\sqrt{r}} \Big\|\Big(\sum_{s=1}^r |k_s|^{-q} |f_{j_s}(\cdot)| ^q\Big)^{\frac{1}{q}}
   \Big\|_{L^{\infty}_x},\\
   \|e^{t\Delta}u_0\|_{\dot{F}^{-1,q}_{\infty}}
   & \lesssim\frac{Q}{\sqrt{r}} \Big\|\Big(\sum_{s=1}^r |k_s|^{-q} |e^{t\Delta}f_{j_s}(\cdot)| ^q\Big)^{\frac{1}{q}}
   \Big\|_{L^{\infty}_x}.
   \end{aligned}
   \label{NS:3.5}
  \end{align}
 Applying \eqref{NS:2.5} to \eqref{NS:2.4}, recalling the properties of $\psi(x)$ and definitions of ${\Psi}_j$ ($j=1,2,3,4$) in
  (H), we have the following point-wise estimates
 \begin{align*}
  \begin{aligned}
  &|\cos(k_sx){\Psi}_1(x)|+|\cos(k_s'x){\Psi}_3(x)|\lesssim (M(\nabla\psi))(x),\\
  &|\sin(k_sx){\Psi}_2(x)|+|\cos(k_s'x){\Psi}_4(x)|\lesssim (M\psi)(x),
   \end{aligned}
 \end{align*}
 where $M\psi$ and $M(\nabla\psi)$ are the Hardy-Littlewood maximal functions of $\psi$ and $\nabla\psi$, respectively.
 Hence
 it follows from \eqref{NS:2.4}-\eqref{NS:2.5}, \eqref{NS:3.5} and Hardy-Littlewood theorem \cite[Chapter\,1, p.13]{Stein:1993} that
 \begin{align*}
    \eqref{NS:3.5}
    &\lesssim
      \frac{Q}{\sqrt{r}}\Big(\Big(\sum_{s=1}^r1\Big)^{\frac{1}{q}}
      \|M\psi\|_{L^{\infty}_x}+\Big(\sum_{s=1}^r|k_s|^{-q}\Big)^{\frac{1}{q}}\|M(\nabla\psi)\|_{L^\infty_x}  \Big)
     \nonumber\\
    &\lesssim  \frac{Q}{\sqrt{r}}(r^{\frac{1}{q}}\|\psi\|_{L^\infty_x}+\|\nabla\psi\|_{L^\infty_x}) \nonumber\\
    &\lesssim  \frac{Q}{\sqrt{r}}(r^{\frac{1}{q}}\|\widehat{\psi}\|_{L^1_\xi}+\|\xi\widehat{\psi}\|_{L^1_\xi}).
  \end{align*}
 Thus \eqref{NS:3.4} follows immediately from  $\|\widehat{\psi}\|_{L^1_\xi}=1$ and  $\|\xi\widehat{\psi}\|_{L^1_\xi}\sim 1$.
\end{proof}
\smallskip
 \begin{lemma}\label{NS:lem3.2}
  For any $T>0$, $u_0$  given in \eqref{NS:2.1}, we obtain that
 \begin{align}\nonumber
   \sup_{0<t<T}t^{\frac{1}{2}}\|e^{t\Delta}u_0\|_{L^{\infty}_x}\lesssim
   {Q}\,{r}^{-\frac{1}{2}}.
 \end{align}
 \end{lemma}
 \begin{proof}
 In view of the initial data construction in \eqref{NS:2.1}, by making use of \eqref{NS:2.4} and $\|M\psi\|_{L^\infty_x}+\|M(\nabla\psi)\|_{L^\infty_x}\lesssim 1$
  as well as $\sup_{t>0}\sum_{s=1}^r e^{-ct2^{2j_s}}t^{\frac{1}{2}}2^{j_s}\lesssim1$, we have
  \begin{align*}
   \sup_{0<t<T}\,t^{\frac{1}{2}}\|e^{t\Delta}u_0\|_{L^\infty_x}
   &\lesssim \frac{Q}{\sqrt{r}}\sup_{t>0}\sum_{s=1}^r
      e^{-ct2^{2j_s}}t^{\frac{1}{2}}(2^{j_s}\|M\psi\|_{L^\infty_x}
      +\|M(\nabla\psi)\|_{L^\infty_x})\\
   &\lesssim    {Q}\,{r}^{-\frac{1}{2}}.
  \end{align*}
  Therefore, we obtain the desired estimate.
 \end{proof}

\smallskip
 Next we need to estimates the norm $\|e^{t\Delta}u_0\|_{X_T}$,
 which is defined as follows
 \begin{align*}
   \|e^{t\Delta}u_0\|_{X_T}=\Big(\sup_{x_0\in\mathbb{R}^3, 0<R^2<T}
   \frac{1}{|B_R(x_0)|}\int_0^{R^2}\int_{B_R(x_0)}|e^{t\Delta}u_0(x)|^2dtdx\Big)^{\frac{1}{2}}.
 \end{align*}
In particular,   $$\|e^{t\Delta}u_0\|_{X_T}\le
\|e^{t\Delta}u_0\|_{L^2_TL^{\infty}_x}.$$

 \begin{lemma}\label{NS:lem3.3}
  For any $T>0$, $u_0$ is given in \eqref{NS:2.1}, for any $0\le N_0\le r$ we
  have
 \begin{align}\nonumber
   \|e^{t\Delta}u_0\|_{X_T}\lesssim
   \frac{Q}{\sqrt{r}}(T^{\frac{1}{2}}|k_{r-N_0}| + \sqrt{N_0}).
 \end{align}
 \end{lemma}
 \begin{proof}
 By the construction of  $u_0$, it suffices to
 estimate $\frac{Q}{\sqrt{r}}{\sum_{\ell}}\,e^{t\Delta}f_{\ell}$.
 Using \eqref{NS:2.4}, $\|f_\ell\|_{L^\infty_x}\lesssim
 |k_s|$ for $\ell=k_s$ and $s\in\{1,2,\cdots,r\}$, and $\|f_\ell\|_{L^\infty_x}=0$ for other $\ell$, we have
\begin{align}
  \frac{Q}{\sqrt{r}}\|\sum_{\ell\in\mathbb{Z}} e^{t\Delta}f_\ell\|_{X_T}&
  \lesssim\frac{Q}{\sqrt{r}}\|\sum_{\ell\in\mathbb{Z}} e^{t\Delta}f_\ell\|_{L^2_TL^{\infty}_x}
  \lesssim  \frac{Q}{\sqrt{r}}\big\| \sum_{s=1}^r e^{-ct2^{2j_s}} 2^{j_s} \big\|_{L^2_T}.\label{NS:3.7}
  \end{align}
 It follows from \eqref{NS:2.12} and
$|k_s|=2^{j_s}$ that for any $N_0\in\![1,r]\cap\, \mathbb{N}$,
 \begin{align}
  \eqref{NS:3.7}
  & \lesssim \frac{Q}{\sqrt{r}}\Big\|\Big(\sum_{s=1}^r |k_s|^{2}e^{-\frac{c}{2}t|k_s|^{2}}\Big)^{\frac{1}{2}}\Big\|_{L^2_T}
    \lesssim \frac{Q}{\sqrt{r}}\Big( \sum_{s=1}^r  \int_0^T\, |k_s|^{2}e^{-\frac{c}{2}t|k_s|^{2}}dt\Big)^{\frac{1}{2}}  \nonumber\\
  & \lesssim \frac{Q}{\sqrt{r}}\Big(\sum_{s=1}^{r-N_0} T|k_s|^{2}\!+\!\!\sum_{s=r-N_0+1}^r\!\!1\Big)^{\frac{1}{2}}
    \lesssim \frac{Q}{\sqrt{r}}({T}^{\frac{1}{2}}|k_{r-N_0}|\!+\!\sqrt{N_0}),\label{NS:3.8}
  \end{align}
 where we used

 $\;\;\;\;\;\;\;\;\;\;\text{$ \int_0^T\, |k_s|^{2}e^{-\frac{c}{2}t|k_s|^{2}}dt \lesssim
\min\{1, \,T|k_s|^{2}\}$
 and $\displaystyle{\sum_{s=1}^{r-N_0}}|k_s|^{2}\!\lesssim\!|k_{r-N_0}|^{2}$.} $ 
\end{proof}

 {\rm By checking the estimate \eqref{NS:3.8} for the case $N_0=0$,  we know that $\|e^{t\Delta}u_0\|_{X_T}\rightarrow 0$ as
  $T\rightarrow 0$. Similarly by checking the estimate \eqref{NS:3.8} for the case $N_0=r$ again, we observe that the best upper bound of
 $\|e^{t\Delta}\|_{X_T}$ is actually $cQ$, which is not good enough to bound the remainder
 $y$ (defined in \eqref{NS:3.1}). 
 Therefore, we need to analyze their contributions by using the
 time-step-division method introduced  in \cite{BP08}.

 Let
 \begin{align}\label{NS:3.9}
 |k_r|^{-2}=T_0 < T_1 < T_2 < \cdots < T_\beta=|k_0|^{-2},
 \end{align}
  where $\beta = Q^3$, $T_\alpha = |k_{r_\alpha}|^{-2}$, $r_\alpha = r -\alpha Q^{-3}r$ and
  $\alpha =0, 1, 2,\cdots,\beta$.

\medskip
\begin{lemma}\label{Lem:3.4}
 Assume that $u_0$ satisfies \eqref{NS:2.1}. Then we have
 \begin{align}
  \label{NS:3.10}
  \Big\|\,e^{t\Delta}u_0\chi_{_{[T_\alpha,T_{\alpha+1}]}}(t)
  \Big\|_{X_{T_{\alpha+1}}}\lesssim \frac{Q}{\sqrt{r}}(1+\sqrt{r Q^{-3}}).
 \end{align}
\end{lemma}
\begin{proof}
  By the construction of initial data $u_0$, we have
  $$
   e^{t\Delta}u_0\chi_{_{[T_\alpha,T_{\alpha+1}]}}(t)
    =
   \frac{Q}{\sqrt{r}}\sum_{\ell\in\mathbb{Z}} \,\chi_{_{[T_\alpha,T_{\alpha+1}]}}(t)\,e^{t\Delta}f_\ell .
  $$
  Similar to \eqref{NS:3.7} and \eqref{NS:3.8}, we get
  \begin{align}\label{NS:3.11}
   \Big\|
     \sum_{\ell\in\mathbb{Z}} \, \chi_{_{[T_\alpha,T_{\alpha+1}]}}(t)\,{e}^{t\Delta} f_\ell
   \Big\|_{X_{T_\alpha+1}}
   \lesssim
    \Big\|
      \chi_{_{[T_\alpha,T_{\alpha+1}]}}(t)\sum_{s=1}^r{e}^{-ct|k_s|^2}|k_s|
    \Big\|_{L^2_{T_{\alpha+1}}}.
  \end{align}
 Applying \eqref{NS:2.12} to \eqref{NS:3.11}, then using Fubini theorem
 and the following facts
 $$\int_{T_\alpha}^{T_{\alpha+1}}
   {e}^{-\frac{c}{2}t|k_s|^2}|k_s|^2dt\lesssim\min\Big\{T_{\alpha+1}|k_s|^2,\,\;1,\;\, e^{-\frac{c}{2}T_{\alpha}|k_s|^2}  \Big\},$$
 we get
 \begin{align}
   \eqref{NS:3.11}
   &\lesssim
     \Big(\sum_{s=1}^{r}\min\{T_{\alpha+1}|k_s|^2,\,\;1,\;\, e^{-\frac{c}{2}T_{\alpha}|k_s|^2}\}  \Big)^{\frac{1}{2}}\lesssim 1+\sqrt{r
     Q^{-3}},\nonumber
 \end{align}
 which can conclude the desired \eqref{NS:3.10}.
\end{proof}
 The following result is a consequence of Lemma \ref{Lem:3.4} since
 $\sum_{s=1}^r e^{-cT_\beta |k_s|^2}\lesssim1$.
\begin{corollary}\label{Cor:3.5}
 For any $T>T_{\beta}=|k_0|^{-2}=2^{-2m_0}$, we have
 \begin{align}
  \label{NS:3.12}
   \|\,e^{t\Delta}u_0\chi_{_{[T_\beta,T]}}(t)
   \,\|_{X_{T}}\lesssim  {Q}{{r}}^{-\frac{1}{2}}.
   \end{align}
\end{corollary}

\subsection{Estimates for the second approximation terms}
 We start this subsection by making some preliminary calculations. In order
  to study the bilinear form  $u_2=B(e^{t\Delta}u_0,e^{t\Delta}u_0)$,
  we
 first split the {\it second approximation terms} $u_2$ into
 \vspace*{-1ex}
 \begin{align}\label{NS:3.13}
 u_{2}={u_{20}}+u_{21}+u_{22},
 \end{align}
 where
   \vspace*{-3.5ex}
    \begin{align}\label{NS:3.14}
  \left\{\begin{aligned}
  u_{20}&=\frac{Q^2}{r}\sum_{s=1}^r \int_0^t
        e^{(t-\tau)\Delta}\mathbb{P}\diva
    \big(
           {e^{\tau\Delta}}{f_{\!j_s}} \otimes {e^{\tau\Delta}}{f}_{\!j_s}
    \big)d\tau, \\
   u_{21}&=\frac{Q^2}{r}\sum_{s=1}^r\sum_{l=1}^{s-1} \int_{0}^t
    {e^{(t-\tau)\Delta}}\mathbb{P}\diva
   \big(
     {e^{\tau\Delta}}{f_{\!j_s}} \otimes {e^{\tau\Delta}}{f}_{\!j_{l}}
   \big)d\tau,\\
     u_{22}&=\frac{Q^2}{r}\sum_{l=1}^r\sum_{s=1}^{l-1} \int_{0}^t
    {e^{(t-\tau)\Delta}}\mathbb{P}\diva
    \big(
     {e^{\tau\Delta}}{f_{\!j_s}} \otimes {e^{\tau\Delta}}{f}_{\!j_{l}}
    \big)d\tau.
    \end{aligned}\right.
 \end{align}

  \noindent \textbf{3.2.1 Analysis of $u_{20}$.}\; To
obtain the lower bound of $u_{20}$, we need to calculate the
 exactly expressions of $u_{20}$ and figure out which part plays the
 key role.
 From \eqref{NS:2.4} and (H), we observe that
  \begin{align*}
   e^{\tau\Delta}f_{j_s}
    \otimes &\,e^{\tau\Delta}f_{j_s}
    =e^{\tau\Delta}(\cos(k_sx)\Psi_1+\cos(k_s'x)\Psi_3)\otimes e^{\tau\Delta}(\cos(k_sx)\Psi_1+\cos(k_s'x)\Psi_3)\\
    &\;\;\;\;\;+|k_s|e^{\tau\Delta}(\cos(k_sx)\Psi_1+\cos(k_s'x)\Psi_3)\otimes{e}^{\tau\Delta}(\sin(k_s'x)\Psi_4-\sin(k_sx)\Psi_2)\\
    &\;\;\;\;\;+|k_s|e^{\tau\Delta}(\sin(k_s'x)\Psi_4-\sin(k_sx)\Psi_2)\otimes{e}^{\tau\Delta}(\cos(k_sx)\Psi_1+\cos(k_s'x)\Psi_3)\\
    &\;\;\;\;\,+|k_s|^2e^{\tau\Delta}(\sin(k_s'x)\Psi_4-\sin(k_sx)\Psi_2)\otimes
    e^{\tau\Delta}(\sin(k_s'x)\Psi_4-\sin(k_sx)\Psi_2)\\
    &\;\;\;\;:=L_{s1}+L_{s2}+L_{s3}+L_{s4}.
  \end{align*}
 Noticing that $L_{s1}$, $L_{s2}$ and $L_{s3}$ are lower order of $|k_s|^2$, hence it suffices
 to estimate
 \begin{align}\label{NS:3.14B}
  L_{s4}=|k_s|^2\,e^{\tau\Delta}\big(\sin(k_s'x){\Psi}_4\!-\!\sin(k_sx){\Psi}_2\big)\otimes{e}^{\tau\Delta}\big(\sin(k_s'x){\Psi}_4\!-\!\sin(k_sx){\Psi}_2\big),
  \end{align}
 {\it where we remark that $L_{s4}$ plays the key role in obtaining the best lower
 bound of $u_{20}$.}
 To be more precisely, by plugging
 $$\sin(k_s'x){\Psi}_4=\frac{e^{ik_s'x}{\Psi}_4-e^{-ik_s'x}{\Psi}_4}{2i}\text{\; and
 \;}\sin(k_sx){\Psi}_2=\frac{e^{ik_sx}{\Psi}_2-e^{-ik_sx}{\Psi}_2}{2i}$$
into \eqref{NS:3.14B}, we can rewrite $L_{s4}$ into the following
four parts:
 \begin{align}\nonumber
  \;\,\,
  \left\{\begin{aligned}
  &J_{s1}=e^{\tau\Delta}(e^{ik_s'x}{\Psi}_4+e^{-ik_sx}{\Psi}_2)\otimes {e}^{\tau\Delta}(e^{ik_s'x}{\Psi}_4+e^{-ik_sx}{\Psi}_2), \\
  &J_{s2}=e^{\tau\Delta}(e^{-ik_s'x}{\Psi}_4+e^{ik_sx}{\Psi}_2)\otimes {e}^{\tau\Delta}(e^{-ik_s'x}{\Psi}_4+e^{ik_sx}{\Psi}_2),\\
  &J_{s3}=e^{\tau\Delta}(e^{ik_s'x}{\Psi}_4+e^{-ik_sx}{\Psi}_2)\otimes {e}^{\tau\Delta}(e^{-ik_s'x}{\Psi}_4+e^{ik_sx}{\Psi}_2), \\
  &J_{s4}=e^{\tau\Delta}(e^{-ik_s'x}{\Psi}_4+e^{ik_sx}{\Psi}_2)\otimes {e}^{\tau\Delta}(e^{ik_s'x}{\Psi}_4+e^{-ik_sx}{\Psi}_2)
  \end{aligned}\right.
\end{align}
such that \vspace*{-2ex}
\begin{align}\label{NS:3.14C}
L_{s4}=J_{s1}+J_{s2}+J_{s3}+J_{s4}.\;
\end{align}
 Correspondingly, we can write that
 \begin{align}\label{NS:3.14D}
 \left\{\begin{aligned}
 &u_{200}=\frac{Q^2}{4r}\sum_{s=1}^r |k_s|^2\int_0^te^{(t-\tau)\Delta}\mathbb{P}\diva
          (J_{s3}+J_{s4})d\tau,\\
 &u_{201}=-\frac{Q^2}{4r}\sum_{s=1}^r |k_s|^2\!\int_0^t\!e^{(t-\tau)\Delta}\mathbb{P}\diva
          (J_{s1}+J_{s2})d\tau,\\
 &u_{202}=\frac{Q^2}{r}\sum_{s=1}^r o(|k_s|^2)\int_0^te^{(t-\tau)\Delta}\mathbb{P}\diva
          (L_{1s}+L_{2s}+L_{3s})d\tau.
 \end{aligned}\right.
 \end{align}

 In what follows, we will spend a lot of effort to deal with $u_{200}$ and get the desired lower bound.
  For any $x_0\in\mathbb{R}^3$,  recalling the definition of ${\Psi}_2$, ${\Psi}_4$ in (H) and
 denoting  $\nu:=k_s+k_s'=(2^3,0,0)$, we have
 \begin{align}\label{NS:3.14E}
  \!\!\!u_{200}(&x_0,t)=C_n\frac{Q^2}{r}\sum_{s=1}^r
  \int_{\mathbb{R}^3\times\mathbb{R}^3}\int_0^t|k_s|^2
  e^{ix_0\xi-(t-\tau)|\xi|^2-\tau(|\xi-\eta|^2+|\eta|^2)}\widehat{\mathbb{P}}\widehat{\diva}
  \nonumber\\
  &\;\times\Big(
  \widehat{\psi}(\xi\!-\!\eta\!-\!k_s')\widehat{\psi}(\eta\!+\!k_s')A_{s}\!+\!\widehat{\psi}(\xi\!-\!\eta\!-\!k_s)\widehat{\psi}(\eta\!+k_s\!)B \nonumber\\
  &\;+
  \widehat{\psi}(\xi\!-\!\eta\!-\!k_s')\widehat{\psi}(\eta\!-\!k_s)D_{s}\!+\!\widehat{\psi}(\xi\!-\!\eta\!+\!k_s')\widehat{\psi}(\eta\!+\!k_s)D_{s}
  \Big) d\tau d\eta d\xi,
 \end{align}
 where $C_n$ is a positive constant depending only on dimension $n$, $\widehat{\diva} =i\xi\,\cdot$,
 $\widehat{\mathbb{P}} $ is a real-valued vector function whose
 $jl$-${\rm th}$ component is
 $\delta_{jl}-\frac{\xi_j\xi_l}{|\xi|^2}$, and 
 \begin{align}\label{NS:3.14F'}
 A_{s}\!=\!\left(\!\!\!
   \begin{array}{ccc}
     2                 & \frac{2^4}{|k_s|}       & 0 \\
     \frac{2^4}{|k_s|} & \frac{2^{7}}{|k_s|^{2}} & 0 \\
     0                 & 0                       & 0 \\
   \end{array}\!
 \right),\;
  B\!=\!\left(
   \begin{array}{ccc}
     0   & 0 & 0 \\
     0   & 0 & 0 \\
     0   & 0 & 2 \\
   \end{array}
 \right),\;
  D_{s}\!=\!\left(\!
   \begin{array}{ccc}
     0                 & 0                  & 1                 \\
     0                 & 0                  & \frac{2^3}{|k_s|} \\
     1                 & \frac{2^3}{|k_s|}  & 0                 \\
   \end{array}\!\!
 \right).
 \end{align}

%

 Next we prove the lower bound of $u_{200}$ in critical space
 $\dot{F}^{-1,q>2}_{\infty}$ and the upper bound of $u_{200}$ in $BMO^{-1}$. Particularly, the lower bound of $u_{200}$ plays
 a {\it crucial} role in the proof of norm inflation. To obtain such bounds,  we will use Fourier analysis methods. Due to the vector-valued nature of velocity field and the divergence free condition, we not only need to explore each of the three components but also need to analyze the action of Leray projection operator
 $\mathbb{P}$. Furthermore, we need to figure out which of the three components is the largest one that produces norm inflation.

 \begin{lemma}{\rm (Lower/upper bound)}\label{Lem:3.11}
 For ${|k_1|^{-2}}\ll T\ll 1$ and $1\le q\le \infty$, we
 have
 \begin{align}\label{DLY-4.34}
 &\|u_{200}(T)\|_{\dot{F}^{-1,q}_{\infty}}\gtrsim{Q}^{2},\\
 &\sup_{0<t<T}\|u_{200}\|_{L^\infty_x}+\|u_{200}\|_{X_T}\lesssim
 T^{\frac{1}{2}}Q^2.\label{DLY-4.35}
 \end{align}
 \end{lemma}
\begin{proof}
 We divide the
 proof of \eqref{DLY-4.34} and \eqref{DLY-4.35} into the following three steps.

 \noindent {\bf{Step 1}.}  Recall that $\nu:=k_s+k_s'=(2^3,0,0)$. From assumptions (H) and \eqref{NS:3.14C}--\eqref{NS:3.14E} we have
 \begin{align}\label{DLY-4.36}
 \textrm{supp}\,\widehat{u}_{200}\!&\subset
  {B}_{\frac{1}{2}}(\nu)
 \cup{B}_{\frac{1}{2}}(-\nu\,)\cup B_{\frac{1}{2}}(0)\subset
 B_{9}(0).
 \end{align}
 Hence for any $t>0$, we have $u_{200}(\cdot,t)\in C^2(\mathbb{R}^3)$. Furthermore, by
 $\dot{F}^{-1,q}_{\infty}\hookrightarrow\dot{B}^{-1,\infty}_{\infty}$,
  \begin{align}
     \|u_{200}\|_{\dot{F}^{-1,q}_{\infty}}\gtrsim\|u_{200}\|_{\dot{B}^{-1,\infty}_{\infty}}
     \gtrsim \|u_{200}\|_{L^\infty_x}\label{4.36-A}
  \end{align}
follows from
 $$\|u_{200}\|_{L^\infty_x}\le \sum_{j\le 3}2^j 2^{-j}\|\Delta_j u_{200}\|_{L^\infty_x}\lesssim
 \sup_{j\in\mathbb{Z}}2^{-j}\|\Delta_ju_{200}\|_{L^\infty_x}=\|u_{200}\|_{\dot{B}^{-1,\infty}_{\infty}}.$$
 We refer readers to \cite[Chapter 5]{LER34} to see more
 information about the equivalence of the definition of Besov
 spaces.

 \noindent \textbf{Step 2}. Considering the arguments in Step 1, it suffices to prove that 
  \begin{align}\label{4.36-B}
    u_{200} (x_0,t)\gtrsim Q^2
  \end{align}
  at some point $x_0$, for instance, here we chose
  $x_0=(\frac{\pi}{2^4},0,0)$.
  Once we prove \eqref{4.36-B}, then combining $u_{200}\in C^2(\mathbb{R}^3)$ with
 \eqref{4.36-A}, we obtain that
  $$\|u_{200}\|_{\dot{F}^{-\alpha,r}_{q_\alpha}}\gtrsim \| u_{200}\|_{L^\infty_x}\gtrsim Q^2,$$
  which is the desired \eqref{DLY-4.34}.

  To prove \eqref{4.36-B}, we first recall from \eqref{NS:3.14} that $u_{20}$ is real-valued which shows that
  the imaginary parts of $u_{200}$, $u_{201}$ and $u_{202}$
  cancels. Therefore, it suffices to bound the real part of
  $u_{200}(x_0,t)$ with $x_0=(-\frac{\pi}{2^4},0,0)$. By \eqref{NS:3.14E}, we set
 \begin{align}\label{NS:3.14EE}
  \Gamma_s=&\,C_n|k_s|^2
  \int_{\mathbb{R}^3\times\mathbb{R}^3}\int_0^t
  (\sin\frac{\pi\xi_1}{2^4})
  e^{-(t-\tau)|\xi|^2-\tau(|\xi-\eta|^2+|\eta|^2)}\widehat{\mathbb{P}}(\xi)
  \nonumber\\
  &\xi\cdot\Big(
  \widehat{\psi}(\xi\!-\!\eta\!-\!k_s')\widehat{\psi}(\eta\!+\!k_s')A_{s}\!+\!\widehat{\psi}(\xi\!-\!\eta\!-\!k_s)\widehat{\psi}(\eta\!+k_s\!)B \nonumber\\
  &+
  \widehat{\psi}(\xi\!-\!\eta\!-\!k_s')\widehat{\psi}(\eta\!-\!k_s)D_{s}\!+\!\widehat{\psi}(\xi\!-\!\eta\!+\!k_s')\widehat{\psi}(\eta\!+\!k_s)D_{s}
  \Big) d\tau d\eta d\xi.
 \end{align}
  It is clear that
  $$({\rm Re\,}u_{200})(x_0,t)=\frac{Q^2}{4r}\sum_{s=1}^r \Gamma_s.$$
 It is also clear that the last two terms in \eqref{NS:3.14EE} are identical
 since substituting $(\xi,\eta)$ by $(-\xi,-\eta)$ yields the same
 results. Furthermore, $A_{s}\rightarrow
 A$
 and $D_{s}\rightarrow D$ as $m_0\rightarrow\infty$, where
 \begin{align}\label{NS:3.14F}
 A\!=\!\left(\!
   \begin{array}{ccc}
     2                 & 0       & 0 \\
     0                 & 0       & 0 \\
     0                 & 0                       & 0 \\
   \end{array}\!
 \right),\;
  B\!=\!\left(
   \begin{array}{ccc}
     0   & 0 & 0 \\
     0   & 0 & 0 \\
     0   & 0 & 2 \\
   \end{array}
 \right),\;
  D\!=\!\left(\!
   \begin{array}{ccc}
     0                 & 0                  & 1                 \\
     0                 & 0                  & 0 \\
     1                 & 0  & 0                 \\
   \end{array}\!
 \right).
 \end{align}

 In order to compute the third component of $\Gamma_s$, for any $\frac{1}{|k_1|^{2}}\ll t\ll 1$, recalling that $\nu=(2^3,0,0)$, then by changing variables, we
 obtain that as $m_0\rightarrow\infty$,
 \begin{align}
 &\int_{\mathbb{R}^3\times\mathbb{R}^3}\int_0^t
  \frac{|k_s|^2 \sin\frac{\pi\xi_1}{2^4}\,
  \widehat{\mathbb{P}}(\xi)\,\xi\!\cdot\!A_s}{e^{t|\xi|^2+\tau(|\xi-\eta|^2+|\eta|^2-|\xi|^2)}}
  \widehat{\psi}(\xi\!-\!\eta\!-\!k_s')\widehat{\psi}(\eta\!+\!k_s')d\tau d\eta d\xi\nonumber\\
  &\;\;\;\;\longrightarrow \int_{\mathbb{R}^3\times\mathbb{R}^3}\frac{\sin\frac{\pi\xi_1}{2^4} }{ e^{t|\xi|^2}} \frac{\xi_3\xi_1^2}{|\xi|^2}\widehat{\psi}(\xi\!-\!\eta)
  \widehat{\psi}(\eta)d\xi d\eta\,; 
  \\
  &\int_{\mathbb{R}^3\times\mathbb{R}^3}\int_0^t
 \frac{|k_s|^2 \sin\frac{\pi\xi_1}{2^4}\,
   \widehat{\mathbb{P}}(\xi)\,\xi\!\cdot\!B}{e^{t|\xi|^2+\tau(|\xi-\eta|^2+|\eta|^2-|\xi|^2)}}
  \widehat{\psi}(\xi\!-\!\eta\!-\!k_s)\widehat{\psi}(\eta\!+\!k_s)d\tau d\eta d\xi\nonumber\\
  &\;\;\;\;\longrightarrow \int_{\mathbb{R}^3\times\mathbb{R}^3}\frac{ \sin\frac{\pi\xi_1}{2^4} } { e^{t|\xi|^2} } \frac{\xi_3(-\xi_1^2\!-\xi_2^2)}{|\xi|^2}\widehat{\psi}(\xi\!-\!\eta)
  \widehat{\psi}(\eta)d\xi d\eta\,;\\
  &\int_{\mathbb{R}^3\times\mathbb{R}^3}\int_0^t
 \frac{|k_s|^2 \sin\frac{\pi\xi_1}{2^4}\,
  \widehat{\mathbb{P}}(\xi)\,\xi\!\cdot\!2D_s}{e^{t|\xi|^2+\tau(|\xi-\eta|^2+|\eta|^2-|\xi|^2)}}
  \widehat{\psi}(\xi\!-\!\eta\!-\!k_s')\widehat{\psi}(\eta\!-\!k_s)d\tau d\eta d\xi\nonumber\\
  &\;\;\;\;\longrightarrow
  \!\int_{\mathbb{R}^3\times\mathbb{R}^3}\!\frac{\sin\frac{\pi(\xi_1\!+\!2^3)}{2^4}}{e^{t|\xi+\nu|^2}}
  \frac{ (\xi_1\!\!+\!2^3)((\xi_1\!+\!2^3)^2\!+\!\xi_2^2\!-\!\xi_3^2)}{|\xi\!+\!\nu|^2}\widehat{\psi}(\xi\!-\!\eta)\widehat{\psi}(\eta)d\xi
  d\eta.
 \end{align}
 Plugging the above three limits and $\|\widehat{\psi}\|_{L^1_\xi}=1$ as well as ${\rm supp}\widehat{\psi}\subset B_{\frac{1}{4}}(0)$ into \eqref{NS:3.14EE},
  for any $2^{-2m_0}\ll t\le 2^{-6}$ and any $k_s$, if $m_0$ is large enough,
  then we get
  $\Gamma_s\sim 1.$
 As a consequence, 
 \begin{align*}
  \|u_{200}(\cdot,t)\|_{L^\infty_x}\gtrsim\|({\rm Re}\,u_{200})(\cdot,t)\|_{L^\infty_x}\gtrsim \frac{Q^2}{r}\sum_{s=1}^r
  \Gamma_s
  \gtrsim {Q^2}.
 \end{align*}
  \textbf{Step 3}. It remains to prove $\|u_{200}\|_{X_T}\lesssim
  Q^2T^{\frac{1}{2}}$ and $\sup_{0<t<T}t^{\frac{1}{2}}\|u_{2,0}\|_{L^\infty_x}$. Using Hausdorff--Young's inequality, we have
 \begin{align*}
 \|u_{200}\|_{{L}^{2}_TL^{\infty}_x}+\sup_{0<t<T}t^{\frac{1}{2}}\|u_{200}\|_{L^\infty_x}\lesssim
 \|\widehat{u}_{200}\|_{L^2_TL^{1}_\xi}+T^{\frac{1}{2}}\|\widehat{u}_{200}\|_{L^1_\xi}.
 \end{align*}
 By checking the proof of Step 2, it is easy to show $\sup_{t>0}|\widehat{u}_{2,0}(\xi,t)|\lesssim Q^2
 (\widehat{\psi}\ast\widehat{\psi})(\xi)$. Consequently,  applying
 Young's inequality to $|\widehat{u}_{2,0}|$ yields
  $$ \|\widehat{u}_{2,0}\|_{L^{1}_\xi} \lesssim
   Q^2 \|\widehat{\psi}\|_{L^{1}_\xi}\|\widehat{\psi}\|_{L^1_\xi}
  \lesssim  Q^2 \quad\text{and }\|\widehat{u}_{2,0}\|_{L^2_TL^{1}_\xi}\lesssim T^{\frac{1}{2}}Q^2 .$$
 Hence we finish the whole proof.\end{proof}
 Now we prove the following estimates for $u_{201}$ and $u_{202}$.
  \begin{lemma}\label{Lem:3.9}
 For any $q>2$, $0<t<T\ll 1$ and large enough $|k_s|$ and $r$, we have
\begin{align}\label{NS:3.15}
 &\|u_{201}\|_{L^\infty_T\dot{F}^{-1,q}_{\infty}}
  +\sup_{0<t<T}t^{\frac{1}{2}}\|u_{201}\|_{L^\infty_x}+\|u_{201}\|_{X_T}
  \lesssim
\frac{Q^2}{\sqrt{r}},\\
&\|u_{202}\|_{L^\infty_T\dot{F}^{-1,q}_{\infty}}
  +\sup_{0<t<T}t^{\frac{1}{2}}\|u_{202}\|_{L^\infty_x}+\|u_{202}\|_{X_T}
  \lesssim
o({Q^2}).
\end{align}
\end{lemma}
\begin{proof}
 First deal with the norm $\|u_{201}\|_{L^\infty_T\dot{F}^{-1,q}_{\infty}}$.
 Noticing from \eqref{NS:2.2B} and \eqref{NS:2.3} that
 \begin{align*}
 &{\rm supp}\widehat{I}_{s1}\subset B_{\frac{1}{2}}(2k_s')\cup B_{\frac{1}{2}}(-2k_s)\cup B_{\frac{1}{2}}(k_s'\!-\!k_s),\\
 &{\rm supp}\widehat{I}_{s2}\subset B_{\frac{1}{2}}(-2k_s')\cup B_{\frac{1}{2}}(2k_s)\cup B_{\frac{1}{2}}(k_s\!-\!k_s').
 \end{align*}
 Hence by $L^\infty_x\hookrightarrow BMO$, $
 BMO^{-1}\hookrightarrow
\dot{F}^{-1,q}_\infty$ and isomorphism as well as boundedness of $\mathbb{P}$ in homogeneous Tribel-Lizorkin spaces, for any $0<t<T$ we get
  \begin{align}\label{NS:3.18}
   \|u_{201}(t)\|_{\dot{F}^{-1,q}_{\infty}}
   &\lesssim
   \frac{Q^2}{r}\big\| \sum_{s=1}^{r} |k_s|^2\int_0^t\!
   e^{(t-\tau)\Delta} (I_{s1}+I_{s2}) d\tau \big\|_{BMO}\nonumber\\&\lesssim
   \frac{Q^2}{r}\big\| \sum_{s=1}^{r} |k_s|^2\int_0^t\!
   e^{(t-\tau)\Delta} (I_{s1}+I_{s2}) d\tau \big\|_{L^\infty_x}.
        \end{align}
 For any $0\le \tau\le t$,  by applying \eqref{NS:2.8} of Lemma \ref{NS:lem2.3} to \eqref{NS:3.18} we  have
    \begin{align}\label{DLY-4.24A}
    \big|e^{(t-\tau)\Delta}(I_{s1}+I_{s2})(x)\big|
 &\lesssim  e^{-ct|k_s|^2}(M\theta)(x) \text{ \; where
 $\theta=(M\psi)^2$}.
 \end{align}
 Plugging \eqref{DLY-4.24A}, $\|M\theta\|_{L^\infty_x}\lesssim 1$ and
 $\displaystyle{\sup_{t>0}\sum_{s=1}^r}\,t|k_s|^2e^{-ct|k_s|^2}\lesssim 1$ into  \eqref{NS:3.18}
 we have
    \begin{align}
   \eqref{NS:3.18}
   &\lesssim \frac{Q^2}{r}
   \sum_{s=1}^{r}\,t|k_s|^2e^{-ct|k_s|^2}  \|M\theta\|_{L^\infty_x} \lesssim
   \frac{Q^2}{r}. \label{DLY-4.25G}
    \end{align}
  To estimate $t^{\frac{1}{2}}\|u_{2,1}\|_{L^\infty_x}$, by using \eqref{DLY-4.24A} and
   $\displaystyle{\sup_{t>0}\sum_{s=1}^r}\,
   t^{\frac{3}{2}}|k_s|^3e^{-ct|k_s|^2}\lesssim1$ we get
 \begin{align}\label{DLY-4.26A}
   t^{\frac{1}{2}}\|u_{2,1}\|_{L^\infty_x}\lesssim\frac{Q^2}{r} \sum_{s=1}^r t^{\frac{3}{2}}|k_s|^3e^{-ct|k_s|^2} \|M\theta\|_{ L^{\infty}_x}
   \lesssim\frac{Q^2}{r}.
  \end{align}
 To estimate $\|u_{2,1}\|_{X_T}$, by using \eqref{DLY-4.24A} and \eqref{NS:2.12} with $\mu=1$ and $\ell=0$
 we get
 \begin{align}\label{DLY-4.26}
   \|u_{2,1}\|_{X_T}&\!\lesssim  \|u_{2,1}\|_{L^2_T L^{\infty}_x} \lesssim\frac{Q^2}{r} \big\|\sum_{s=1}^r t|k_s|^3e^{-ct|k_s|^2}\big\|_{L^2_T}
   \big\|M\theta\big\|_{ L^{\infty}_x}\nonumber\\
   &\lesssim\frac{Q^2}{r}\big\|\big(\sum_{s=1}^r |k_s|^2e^{-\frac{c}{2}t|k_s|^2}\big)^{\frac{1}{2}} \big\|_{L^2_T}
   \lesssim \frac{Q^2}{\sqrt{r}}.
  \end{align}
Combining \eqref{NS:3.18} and \eqref{DLY-4.25G}--\eqref{DLY-4.26},
we finish the proof of
 \eqref{NS:3.15}.

 To estimate  $u_{202}$, we recall
 that similar to \eqref{NS:3.14B},
 \begin{align*}
  u_{202}\sim \frac{Q^2}{r}\sum_{s=1}^r o(|k_s|^2)\int_0^t e^{(t-\tau)\Delta}\mathbb{P}\diva (L_{1s}+L_{2s}+L_{3s}) d\tau,\end{align*}
 where  $o(|k_s|^2)$ is lower order term of $|k_s|^2$ satisfying $0\le \frac{o(|k_s|^2)}{|k_s|^2}\lesssim 2^{-m_0}$.
 Similar to $L_{4s}$, for large enough $|k_0|$ and $r$, we can prove the desired estimate.
%
\end{proof}


\noindent\textbf{3.2.2 Analysis of $u_{21}$ and $u_{22}$}
 \begin{lemma}\label{Lem:3.10}
 For any $q>2$ and $i=1,2$, we have
\begin{align}\label{DLY-4.27}
 \|u_{2i}\|_{L^\infty_T\dot{F}^{-1,q}_{\infty}}+\sup_{0<t<T}t^{\frac{1}{2}}\|u_{2i}\|_{L^\infty_x}+\|u_{2i}\|_{X_T}
  \lesssim {Q^2}{{r}}^{-1}.
\end{align}
\end{lemma}
\begin{proof}
 It suffices to estimate $u_{21}$. To estimate
 $\|u_{21}\|_{L^\infty_T\dot{F}^{-1,q}_{\infty}}$, recalling from the
 support of $\widehat{f}_{j_s}$ and $|k_l|=j_l\ne j_s=|k_s|$, we get
   \begin{align*}
   &\supp \mathcal{F}(e^{\tau\Delta}f_{j_s}\otimes e^{\tau\Delta}f_{j_l})\subset
    B_{\frac{1}{2}}(\pm k_s'\pm k_l)\cup B_{\frac{1}{2}}(\pm k_s\pm k_l').
   \end{align*}
 Note that when $s>\ell$, $|\pm k_s'\pm k_l|\!\sim\!|k_s|\!\sim\!|\pm k_s\pm
 k_l'|$, thus the support of $\widehat{u}_{21}$ is far away from origin which ensures that $\mathbb{P}$ is well-defined and no singularity arguments for $\mathbb{P}$ are involved.
 Moreover, $\mathbb{P}$ is a bounded operator in $BMO$. By $L^\infty_x\!\hookrightarrow\!BMO$, $BMO^{-1}\!\hookrightarrow\!\dot{F}^{-1,q}_\infty$, isomorphism,
  Lemma \ref{NS:lem2.3} and (H), we get
  \begin{align}
   \|u_{21}\|_{L^\infty_T\dot{F}^{-1,q}_{\infty}}
   &\lesssim
   \frac{Q^2}{r}\big\|
   \sum_{s=1}^{r}
   \sum_{l=1}^{s-1}
   \int_0^t
    e^{(t-\tau)\Delta}(e^{\tau\Delta}f_{j_s}\otimes
    e^{\tau\Delta}f_{j_l})
    d\tau\big\|_{L^\infty_TBMO}\nonumber\\
   &\lesssim
   \frac{Q^2}{r}\big\|
   \sum_{s=1}^{r}
   \sum_{l=1}^{s-1}
   \int_0^t
    e^{(t-\tau)\Delta}(e^{\tau\Delta}f_{j_s}\otimes
    e^{\tau\Delta}f_{j_l})
    d\tau\big\|_{L^\infty_TL^\infty_x}\nonumber\\
   &\lesssim\frac{Q^2}{r}\sup_{t>0}\sum_{s=1}^r\sum_{l=1}^{s-1}\int_{0}^t e^{-ct|k_s|^2} d\tau \nonumber\\
   &\hskip1cm \big\|M\big(\, (\,|k_s|M\psi+M(\nabla\psi)(\,|k_l|M\psi+M(\nabla\psi)\,\big)\big\|_{L^\infty_x}\nonumber\\
   &\lesssim \frac{Q^2}{r}\sup_{t>0}\sum_{s=1}^r\sum_{l=1}^{s-1} e^{-ct|k_s|^2}t|k_s||k_l|\nonumber\\
   &\lesssim \frac{Q^2}{r},\label{NS:3.25}
        \end{align}
 where in the  fourth inequalities we used the following simple
 fact:
 \begin{align}\label{NS:3.26}
  \big\|M\big((|k_s|M\psi\!+\!M(\nabla\psi)(|k_l|M\psi\!+\!M(\nabla\psi)\big)\big\|_{L^\infty_x}\lesssim
  |k_s||k_l|.
 \end{align}
%
 Next to the norm of $\|u_{21}\|_{X_T}$. By using Lemmas \ref{NS:lem2.3} and \ref{NS:lem2.4}, then we have
 \begin{align}\label{DLY-4.29}
 \|u_{21}\|_{X_T} &\lesssim \|u_{21}\|_{L^2_T L^{\infty}_x}
 \lesssim \frac{Q^2}{r}\Big\| \sum_{s=1}^r\sum_{l=1}^{s-1} t|k_s|^2|k_{l}|e^{-ct|k_s|^2} \Big\|_{L^2_T}\nonumber\\
 &\lesssim\frac{Q^2}{r}\Big\|\Big(\sum_{s=1}^r |k_{s-1}|^2 e^{-\frac{c}{2}t|k_s|^2}\Big)^{\frac{1}{2}} \Big\|_{L^2_T}
    \lesssim \frac{Q^2}{r}.
  \end{align}
  At last, we estimate $t^{\frac{1}{2}}\|u_{21}\|_{L^\infty_x}$. By using Lemmas \ref{NS:lem2.3} and \ref{NS:lem2.4}, we have
 \begin{align}\label{DLY-4.29A}
  \sup_{0<t<T} t^{\frac{1}{2}}\|u_{21}\|_{L^\infty_x}
   &\lesssim\frac{Q^2}{r} \sum_{s=1}^r\sum_{l=1}^{s-1}
   t^\frac{3}{2}|k_s|^2|k_{l}|e^{-ct|k_s|^2}
   \lesssim \frac{Q^2}{r}.
  \end{align}

 Thus combining \eqref{NS:3.25}, \eqref{DLY-4.29} and \eqref{DLY-4.29A}, we finish the
 proof of \eqref{DLY-4.27}.
\end{proof}
\subsection{Estimates of remainder $y$}
 In this subsection, we use iteration arguments to prove the {\it a-priori} estimate for remainder $y$.
 Recall that $y$  satisfy the integral equations
 \eqref{NS:3.2}, i.e.
 $$y=G_0+G_1-G_2$$
  with initial condition $y|_{t=0}=0$, $ G_0  = B(u_{2}, u_1)+B(u_1,u_{2})- B(u_{2},u_{2})$ and
  \begin{align*}
  \begin{aligned}
  G_1 & = B(y,u_2-u_1)+B(u_2-u_1,y), \quad\
 G_2  = B(y,y).
      \end{aligned}
      \end{align*}

 From Lemma \ref{Lem:3.4}, we observe that in order to obtain more accurate decay estimate for $y$,
 it suffices to split $u_1$, $u_2$ and $u_{2}$ into
 two terms, e.g.
 \begin{align*}
 \left\{\begin{aligned}&u_1=u_1\chi_{[0,T_\alpha]}(t) +
 u_1\chi_{[T_\alpha,T_{\alpha+1}]}(t),\\
 &u_2=u_2\chi_{[0,T_\alpha]}(t) +
  u_2\chi_{[T_\alpha,T_{\alpha+1}]}(t),\\
 &y=y\chi_{[0,T_\alpha]}(t) + y\chi_{[T_\alpha,T_{\alpha+1}]}(t),
\end{aligned}\right.
\end{align*}
 Plugging the above decompositions of $u_1$, $u_2$ and $y$ into
 $G_0$, $G_1$ and $G_2$, we have the following iteration rules which play an important role in controlling $y$.
\begin{lemma}\label{Lem:3.12}
 If $y$ solves system \eqref{NS:3.2}, then for any $\alpha=0,1,2,\cdots,Q^3$ and for large
 enough $r$ and $|k_0|$ we have
 \begin{align}\label{DLY:3.71}
 \|y\|_{X_{T_{\alpha+1}}}%
    \lesssim Q^{\alpha+3} (\, {r}^{-\frac{1}{2}}+ {|k_0|^{-1}}\, ).
    \end{align}
 Moreover, for any $T>|k_0|^{-2}$, we have
 \begin{align}\label{DLY:3.72}
 \|y\|_{X_{T}}
     \lesssim Q^{3} ( {r}^{-\frac{1}{2}} + T^{\frac{1}{2}} )+Q^{Q^3+3} (\, {r}^{-\frac{1}{2}} +  {|k_0|^{-1}}\, ).
    \end{align}
 \end{lemma}

 \begin{proof}
 Applying  \ref{NS:1.8} to \eqref{NS:3.2}, we
 have the following bilinear estimates:
\begin{align}
   \|y\|_{X_{T_{\alpha+1}}}
   &\lesssim  \|u_2\|_{X_{T_{\alpha+1}}}(\|u_1\|_{X_{T_{\alpha+1}}} \!+\! \|u_2\|_{X_{T_{\alpha+1}}}) + (\|u_1\|_{X_{T_{\alpha+1}}}\!\!\!+\!\|u_2\|_{X_{T_{\alpha+1}}})\|y\|_{X_{T_\alpha}}\nonumber\\
   &\;\;\;+(\|u_1\|_{X_{[T_\alpha,T_{\alpha+1}]}}\!+\!\|u_2\|_{X_{[T_\alpha,T_{\alpha+1}]}})\|y\|_{X_{T_{\alpha+1}}}\nonumber\\
   &\;\;\;+\|y\|_{X_{T_{\alpha+1}}}^2,\label{DLY:3.73}
     \end{align}
 where in the second inequality we used
 $$\|y\|_{X_{T_{\alpha}}}\!+\|y\|_{X_{[T_\alpha,T_{\alpha+1}]}}
 \lesssim
 \|y\|_{X_{T_{\alpha+1}}}\;\text{ and }\;
 \|u_2\|_{X_{T_\alpha}}\!+\|u_2\|_{X_{[T_{\alpha},T_{\alpha+1}]}} \!\lesssim\! \|u_2\|_{X_{T_{\alpha+1}}}.$$
 Recalling that for any $1\le\alpha\le\beta$, $T_\alpha\le T_\beta$. Then from Lemmas \ref{NS:lem3.2}--\ref{Lem:3.11}, we get
 \begin{align}\label{DLY:3.74}
 \|u_2\|_{X_{T_{\alpha+1}}}\!\lesssim\!
 {Q^2}{{r}}^{-\frac{1}{2}} \!+ Q^2T^{\frac{1}{2}}_\beta,\;\;
 \|u_1\|_{X_{T_\alpha}} \!\lesssim\! Q,\;\; \|u_1\|_{X_{[T_\alpha,T_{\alpha+1}]}}
 \!\lesssim\! Q^{-\frac{1}{2}}.
 \end{align}
 Plugging \eqref{DLY:3.74} in \eqref{DLY:3.73}, and assuming that $r>Q^{10}$, $T_\beta=|k_0|^{-2}<Q^{-5}$, we have
 \begin{align}\label{DLY:3.75}
 \|y\|_{X_{T_{\alpha+1}}}
  &\lesssim Q^3(\, {{r}^{-\frac{1}{2}}} + {|k_0|^{-1}}\,)
             +Q\|y\|_{X_{T_\alpha}} + Q^{-\frac{1}{2}}\|y\|_{X_{T_{\alpha+1}}} + \|y\|_{X_{T_{\alpha+1}}}^2.
 \end{align}
 Similarly, when $T>T_\beta$,
 by splitting $[0,T]$ into $[0,T_\beta]$ and $[T_\beta,T]$, then
 using
  Corollary \ref{Cor:3.5} and \eqref{DLY:3.73}--\eqref{DLY:3.75}, we get $\|u_1\|_{X_{[T_\beta,T]}}\lesssim
 Qr^{-\frac{1}{2}}$ and
 \begin{align}\label{DLY:3.76}
 \|y\|_{X_{T}}
 &\lesssim Q^3(\,{{r}^{-\frac{1}{2}}} +T^{\frac{1}{2}})
             +Q\|y\|_{X_{T_\beta}} + Q^{-\frac{1}{2}}\|y\|_{X_{T}} + \|y\|_{X_{T}}^2.
 \end{align}
 Lemma \ref{NS:lem3.2} ensures that $\|y\|_{X_{T_0}}$ can be small
 since $T_0=|k_r|^{-2}$ and
 $$\|u_1\|_{X_{T_0}}\lesssim
 \frac{Q}{\sqrt{r}}T^\frac{1}{2}_0|k_r|\lesssim \frac{Q}{\sqrt{r}}$$ which can be arbitrarily small if $r$ is large
 enough. Thus iteration argument can be applied to \eqref{DLY:3.75}--\eqref{DLY:3.76}. Iterating
 \eqref{DLY:3.75} and \eqref{DLY:3.76} give the desired results.
 \end{proof}

 Making use of \eqref{NS:1.8} and Lemma \ref{Lem:3.12}, we obtain the
 following estimate.
\begin{corollary}\label{Cor:3.13}
 For any  $q>2$, sufficiently large $r$ and $|k_0|$ such that $r\gg Q^{2Q^3+4}$, $|k_0|\gg Q^{-Q^3-2}$
 and $|k_0|^{-2}<T\ll Q^{-2}$, we have
 \begin{align}\label{DLY:3.77}
   & \|y(T)\|_{\dot{F}^{-1,q}_{\infty}}
\ll Q^2.
 \end{align}
\end{corollary}
\begin{proof}
 From \eqref{NS:3.2}, we notice that $y(T)=G_0(T)+G_1(T)-G_2(T)$
 and $G_i(T)$ are several bilinear terms. By   \eqref{NS:1.8} and \eqref{NS:1.9}, we obtain that
 \begin{align}
 \|y(T)\|_{\dot{F}^{-1,q}_{\infty}}
   &\lesssim \|y(T)\|_{\dot{F}^{-\alpha,2}_{\infty}}\lesssim   \|y\|_{L^\infty_T\dot{F}^{-\alpha,2}_{\infty}}\nonumber\\
   &\lesssim
   \|u_2\|_{X_T}(\|u_1\|_{X_T}+\|u_2\|_{X_T})+\|y\|_{X_T}(\|u_1\|_{X_T}+\|u_2\|_{X_T})+\|y\|_{X_T}^2.\nonumber
 \end{align}
  Applying Lemmas \ref{NS:lem3.2},
 \ref{Lem:3.9}--\ref{Lem:3.12} to the above inequality, we have
 \begin{align}
  \|y(T)\|_{\dot{F}^{-1,q}_{\infty}}
  &\lesssim( {Q^2}{r}^{-\frac{1}{2}} +  Q^2T^{\frac{1}{2}})(Q +  {Q^2}{r}^{-\frac{1}{2}} + Q^2T^{\frac{1}{2}})
  \nonumber\\
  &\;\;\;\;+(Q +   {Q^2}{r}^{-\frac{1}{2}} + Q^2T^{\frac{1}{2}})
  \Big(Q^{3}( {r}^{-\frac{1}{2}} + T^{\frac{1}{2}})+Q^{Q^3+3}( {r}^{-\frac{1}{2}} + {|k_0|^{-1}})\Big)\nonumber\\
  &\;\;\;\;+\Big(Q^{3}( {r}^{-\frac{1}{2}} +\!T^{\frac{1}{2}})+Q^{Q^3+3}(\, {r}^{-\frac{1}{2}} + {|k_0|^{-1}})\Big)^2
  \nonumber\\
  &\ll Q^2.\nonumber \end{align}
 Hence we prove the desired result.
  \end{proof}

 \subsection{Proof of Theorem
 \ref{Thm:1.3}}
  In this subsection, combining the results proved in Subsections 3.1--3.4, we are ready to prove the norm inflation of the Navier-Stokes equations.

\medskip
 \noindent\textit{Proof of Theorem
 \ref{Thm:1.3}}. Combining the equalities \eqref{NS:3.1} and \eqref{NS:3.13}, the estimates \eqref{NS:3.4},
 \eqref{NS:3.15}, \eqref{DLY-4.27},  \eqref{DLY-4.34} and \eqref{DLY:3.77}, we have
 \begin{align}
   \|u(T)\|_{\dot{F}^{-1,q}_{\infty}}&\ge
   \|u_{200}(T)\|_{\dot{F}^{-1,q}_{\infty}}
   - \|u_1(T)\|_{\dot{F}^{-1,q}_{\infty}}\!\nonumber\\
   &\;\;\;-\!\sum_{\ell=1}^2(\|u_{20\ell}(T)\|_{\dot{F}^{-1,q}_{\infty}}+\|u_{2\ell}(T)\|_{\dot{F}^{-1,q}_{\infty}})
   \!-\!\|y(T)\|_{\dot{F}^{-1,q}_{\infty}}\nonumber\\
   &\gtrsim Q^2\Big(1-{Q}^{-1}r^{\frac{1}{q}-\frac{1}{2}}-r^{-\frac{1}{2}}-o(1)\Big)
    \gtrsim Q^2,\nonumber
 \end{align}
 where $0<o(1)\ll\frac{1}{2}$, $r\!\gg\!Q^{2Q^3+4}$ and $|k_0|^{-2}\!<\!T\!\ll \!Q^{-2}$.
  Hence we finish the proof.

 \vspace*{2ex}
{{\small \hspace*{4.3cm}Chao Deng,}\\
 {\small \hspace*{5cm}Department of Mathematics,}\\
 {\small \hspace*{5cm}Jiangsu Normal University, Xuzhou 221116, China;}\\
 {\small \hspace*{5cm}Department of Mathematics, Penn State University.}\\
 {\small \hspace*{5cm}{\it{deng315@yahoo.com.cn}}}\\
\\
{\small \hspace*{5cm}Xiaohua Yao,}\\
 {\small \hspace*{5cm}Department of  Mathematics,}\\
 {\small \hspace*{5cm}Central China Normal University, Wuhan 430079, China.}\\
 {\small \hspace*{5cm}{\it{yaoxiaohua@mail.ccnu.edu.cn}}}\\
}
 \date{}

\end{document}